\documentclass{proc-l}

\newtheorem*{theorem}{Theorem}
\newtheorem*{corollary}{Corollary}

\theoremstyle{definition}
\newtheorem*{definition}{Definition}

\theoremstyle{remark}
\newtheorem*{remark}{Remark}
\newtheorem*{ack}{Acknowledgment}

\numberwithin{equation}{section}
\begin{document}

\title{Separable subgroups have bounded packing}


%    author one information
\author{Wen-yuan Yang}
\address{College of Mathematics and Econometrics, Hunan
University, Changsha, Hunan 410082 People's Republic of China}
\curraddr{U.F.R. de Mathematiques, Universite de Lille 1, 59655
Villeneuve D'Ascq Cedex, France} \email{wyang@math.univ-lille1.fr}
\thanks{The author is supported by the China-funded Postgraduates
Studying Aboard Program for Building Top University. This research
was supported  by National Natural Science Foundational of China
(No. 11071059)}

%    author two information

%    \subjclass is required.
\subjclass[2000]{Primary 20F65, 20F67}
%    The 2010 edition of the Mathematics Subject Classification is
%    now available.  If you are citing a classification from the
%    new scheme, use the following input coding instead.
%\subjclass[2010]{Primary }

\date{}

\dedicatory{}

%    "Communicated by" -- provide editor's name; required.
\commby{Alexander Dranishnikov}

\keywords{bounded packing, separable, polycyclic groups}

%    Abstract is required.
\begin{abstract}
In this note, we prove that separable subgroups have bounded packing
in ambient groups. The notion bounded packing was introduced by
Hruska and Wise and in particular, our result answers positively a
question of theirs, asking whether each subgroup of a virtually
polycyclic group has the bounded packing property.
\end{abstract}

\maketitle

\section{Introduction}
Bounded packing was introduced for a subgroup of a countable group
in Hruska-Wise \cite{HrWi}. Roughly speaking, this property gives a
finite upper bound on the number of left cosets of the subgroup that
are pairwise close in $G$. Precisely,
\begin{definition}
Let $G$ be a countable group with a left invariant proper metric
$d$. A subgroup $H$ has \textit{bounded packing} in $G$ (with
respect to $d$) if for each positive constant $D$, there is a
natural number $N=N(G,H,D)$ such that, for any collection
$\mathcal{C}$ of $N$ left $H$-cosets in $G$, there exist at least
two $H$-cosets $gH,g'H \in \mathcal{C}$ satisfying $d(gH,g'H) > D$.
\end{definition}

\begin{remark}
Bounded packing of a subgroup is independent of the choice of the
left invariant proper metric $d$. Equivalently, bounded packing says
that for each positive constant $D$, every collection of left
$H$-cosets in $G$ with pairwise distance at most $D$ has a uniform
bound $N=N(G,H,D)$ on their cardinality.
\end{remark}

This note aims to give a proof of the following.
\begin{theorem}\label{thm:separable}
If $H$ is a separable subgroup of a countable group $G$, then $H$
has bounded packing in $G$.
\end{theorem}

A subgroup $H$ of a group $G$ is \textit{separable} if $H$ is an
intersection of finite index subgroups of $G$. A group is called
\textit{subgroup separable} or \textit{LERF} if every finitely
generated subgroup is separable. For example, Hall showed that free
groups are LERF in \cite{Hall}. It follows from a theorem of Mal'cev
\cite{Mal} that polycyclic (and in particular finitely generated
nilpotent) groups are LERF. A group is called \textit{slender} if
every subgroup is finitely generated. Polycyclic groups are also
slender by a result of Hirsch \cite{Hir}.  Therefore, we have the
following corollary, which gives a positive answer to
\cite[Conjecture 2.14]{HrWi}.
\begin{corollary}
Let $P$ be virtually polycyclic. Then each subgroup of $P$ has
bounded packing in $P$.
\end{corollary}
\begin{remark}
In \cite{Sah}, Jordan Sahattchieve obtained a special case of this Corollary using different methods: any subgroup of (Hirsch) length 1 of a polycyclic group has bounded packing.  
\end{remark}

\ack The author would like to sincerely thank Prof. Leonid
Potyagailo for his comments and interest in this work.

\section{Proof of the Theorem}

We define the norm $|g|_d$ of an element $g \in G$ as the distance
$d(1,g)$.

\begin{proof}[Proof of the Theorem]
By the definition of bounded packing, it suffices to show, for each
positive constant $D$, that there is a uniform bound on the
cardinality of every collection of left $H$-cosets in $G$ with
pairwise distance at most $D$.

Given such a collection $\mathcal{A}$ satisfying $d(gH,g'H) < D$ for
any $gH$, $g'H \in \mathcal{A}$. Without loss of generality, we can
assume $H$ belongs to $\mathcal{A}$, up to a translation of
$\mathcal{A}$ by an appropriate element of $G$. Since $d(H,gH)<D$
for each $gH \in \mathcal{A}$, there exists an element $h$ in $H$
such that $d(1,hgH) < D$. Hence we conclude that the collection
$\mathcal{A} \setminus \{H\}$ lies in the finite union of double
cosets $HgH$ with $|g|_d<D$ and $g \in G \setminus H$.

Since $d$ is a left invariant proper metric on $G$, the set $F = \{g
\in G \setminus H: |g|_d<D\}$ is finite. Since $H$ is separable in
$G$, we can take a finite index subgroup $K$ of $G$ such that $H<K$
and $F \subset G \setminus K$.

We claim that no two different left $H$-cosets of $\mathcal{A}$ lie
in the same left $K$-coset. By way of contradiction, we suppose that
there are two $H$-cosets $g k H$, $g k' H \in \mathcal{A}$ in the
same coset $g K$ such that $d(g k H,g k' H) < D$. By a similar
argument as above, we get that $k^{-1}k'H$ belongs to a double coset
$Hg_0 H$ with $|g_0 |_d<D$. Moreover, we note that $g_0 \in F$.
Since we have $k^{-1}k'H = hg_0 H$ for some $h \in H$, it is easy to
see that $g_0$ belongs to $K$. But by the choice of $K$, we know
that $g_0$ belongs to $G \setminus K$. This is a contradiction. Our
claim is proved.

Since $K$ is of finite index in $G$, the cardinality of each
$\mathcal{A}$ is upper bounded by $[G:K]$. Thus for each $D$, we
have obtained a uniform bound on every $\mathcal{A}$. Hence $H$ has
bounded packing in $G$.
\end{proof}

\bibliographystyle{amsplain}

\begin{thebibliography}{10}

\bibitem{Hall} M. Hall, \textit{Coset representations in free groups.} Trans. Amer. Math. Soc., \textbf{67} (1949), 421--432.
\bibitem{Hir} K. A. Hirsch, \textit{On infinite soluble groups I.} Proc. London Math. Soc., \textbf{s2-44(1)} (1938), 53--60.
\bibitem{HrWi} G. Hruska and D. Wise, \textit{Packing subgroups in relatively hyperbolic groups. } Geom. Topol.,  \textbf{13(4)} (2009), 1945--1988.
\bibitem{Mal} A. I. Mal'cev, \textit{On homomorphisms onto finite groups.} American Mathematical Society Translations, \textbf{119(2)} (1983), 67--79.
\bibitem{Sah} J. Sahattchieve, \textit{On Bounded Packing in Polycyclic Groups.} Preprint,  arXiv:1011.5953, 2010.
 

\end{thebibliography}

\end{document}